\documentclass[11pt]{amsart}

\usepackage{amsfonts}
\usepackage{amscd}
\usepackage{amssymb}
\usepackage{amsmath}
\usepackage{latexsym}
\usepackage[dvips]{graphicx}
\usepackage{rotating}
\usepackage{tabularx}

\newtheorem{prop}{Proposition}[section]
\newtheorem{defn}[prop]{Definition}
\newtheorem{lem}[prop]{Lemma}
\newtheorem{rem}[prop]{Remark}

\newtheorem{theo}[prop]{Theorem}

\def\CQFD{$\Box$}
\newenvironment{preuve}{{\it Proof}.}{\hfill \CQFD \\}

\def\KK{{\mathbb{K}}}
\def\PP{{\mathbb{P}}}
\def\RR{{\mathbb{R}}}

\def\NN{{\mathbb{N}}}
\def\ZZ{{\mathbb{Z}}}

\def\OO{{\mathcal{O}}}
\def\LL{{\mathcal{L}}}

\def\MM{{\mathcal{M}}}

\def\HH{{\mathcal{H}}}
\def\XX{{\mathcal{X}}}
\def\UU{{\mathcal{U}}}
\def\DD{{\mathcal{D}}}
\def\EE{{\mathcal{E}}}
\def\FF{{\mathcal{F}}}

\def\dg{\mathbf{d}}
\def\kg{\mathbf{k}}

\def\Res{\mathrm{Res}}
\def\End{\mathrm{End}}
\def\dim{\mathrm{dim}}
\def\det{\mathrm{det}}
\def\deg{\mathrm{deg}}

\def\Spec{\mathrm{Spec}}
\def\codim{\mathrm{codim}}

\def\rank{\mathrm{rank}}

\def\Hom{\mathcal{H}om}
\def\hom{\mathrm{Hom}}
\def\Im{\mathrm{Im}}
\def\coker{\mathrm{coker}}

\def\rth{r^{\scriptsize{th}}}

\title[Determinantal resultant]{Determinantal resultant}
\date{August 15, 2002 }
\author{Laurent Bus\'e}
\address{ Universit\'e de Nice - Sophia Antipolis,\newline Laboratoire de
  Math\'ematiques, Parc Valrose 06108 Nice Cedex 2.} 
\email{lbuse@unice.fr}

\begin{document}
\maketitle

\begin{abstract} In this paper, 
 a new kind of resultant, called the determinantal resultant, is
 introduced. This operator computes the projection of a determinantal
 variety under suitable hypothesis. As a direct generalization of 
 the resultant of a very ample vector bundle \cite{GKZ94}, it corresponds to a
 necessary and sufficient condition so that a given morphism between
 two vector bundles on a projective variety $X$ has rank lower or equal to
 a given integer in at least one point. First some conditions are given for the existence of such a
 resultant and it is showed how to compute explicitly its degree. Then
 a  result of A. Lascoux \cite{Las78} is used  to obtain it 
 as a 
 determinant of a certain complex. Finally some more detailed 
 results in the particular case 
 where $X$ is a projective space are exposed. 
\end{abstract}

\section{Introduction}
Projection is one of the more used operation in effective algebraic
geometry and more particularly in elimination theory. The Gr\"obner basis
theory is a powerful tool to perform projections in all generality,
however its computation complexity 
can be very high in practice. It is hence useful to develop some
other tools being able to compute such projections, even if they apply
only in particular cases. \emph{Resultants} 
are such tools: they are used to eliminate a set of variables (often
called parameters) of a given polynomial system, providing this
elimination process leads to only one equation (in the
parameters). The most known resultant is the  
Sylvester resultant of two homogeneous bivariate polynomials. Its
generalization to $n$ homogeneous polynomials in $n$ variables has
been stated by F.S. Macaulay in 1902 \cite{Mac02}. Let 
$f_1,\ldots,f_n$ be polynomials in variables $x_1,\ldots,x_n$, their 
resultant is a polynomial in their coefficients  
which vanishes if and only if they have a common homogeneous root (in
the algebraic closure of the ground field). An improvement of this
resultant, taking account of the
monomial supports of the input polynomials, has been exposed in
\cite{KSZ92} and yields a tool even more efficient called the sparse
resultant (see also \cite{Stu93}).  

In the book of
I.M. Gelfand, M.M. Kapranov and   A.V. Zelevinsky \cite{GKZ94}, a more
general resultant, which encapsulates all the previous cited
resultants, is introduced~: the resultant of 
a very ample vector bundle $E$ of rank $n$ on an  irreducible projective 
variety $X$ of dimension $n-1$ over an algebraically closed field
$\KK$. Denoting $V=H^0(X,E)$, this resultant is the divisor of the
projective space $Y=\PP(V)$ obtained as the projection on $Y$ of the
incidence variety 
$$W=\{ (x,f) \in X\times Y \  :  \  f(x)=0 \},$$
 by the canonical map 
$X\times Y\rightarrow Y$. For instance, the Macaulay resultant
corresponds to the case $X=\PP^{n-1}$ and
$E=\oplus_{i=1}^n\OO_{\PP^{n-1}}(d_i)$ with $d_i>0$ for
$i=1,\ldots,n$ (in such a situation where the vector bundle $E$ splits
into $n$ line bundles the resultant is called the \emph{mixed}
resultant in \cite{GKZ94}). The aim of this paper is to 
extend this resultant. Being given $n$ homogeneous polynomials 
$f_1,\ldots,f_n$ in variables $x_1,\ldots,x_n$, the Macaulay resultant
vanishes if and only if 
$f_1,\ldots,f_n$ have a common root on $\PP^{n-1}$,
that is the line matrix $(f_1,\ldots,f_n)$ is of rank 
zero in at least one point of $\PP^{n-1}$. More generally the
resultant of the vector bundle $E$ can be interpreted 
as an operator which traduces a rank default of a morphism 
$$E^* \xrightarrow{ \ \,  f \ \, }  \OO_X,$$ in at least one point of
$X$. Consequently, we can ask for a ``resultant'' condition so that a
morphism between two vector bundles on an irreducible projective 
variety $X$ is of rank lower or equal to a given integer (which was
for instance zero in the previous example). 
The first apparition of this problem in a concrete situation seems to be 
a talk of H. Lombardi. Notice that J.P. Jouanolou has also worked on 
such a problem, apparently in the particular case that we will call
the principal case (unpublished). From an historically point of view
the first use  
of determinantal ideals and Eagon-Northcott complexes (that we will
encounter 
hereafter) in elimination theory was probably the paper \cite{Laz77} of
D. Lazard.     

This paper is organized as follows. In the section \ref{secdefdetres},
 we define 
what we will call a $\rth$-determinantal resultant. Being given 
two vector bundles $E$ and $F$ of respective rank 
$m\geq n$ on an irreducible  projective variety $X$ and an  
integer $0\leq r <n$, the $\rth$-determinantal resultant of
$E$ and $F$ on $X$ is a divisor 
on the vector space of morphisms of $E$ in $F$ which gives a 
necessary and sufficient condition so that such a morphism has rank
lower or equal to $r$ in at least on point of $X$. We give suitable 
hypothesis on the dimension of $X$ and on the vector
bundles $E$ and $F$ such that this determinantal resultant
exists. Then we show that, 
as for the resultant of a very ample vector bundle of \cite{GKZ94}, we can
give the degree (and the multi-degree if $E$ splits) of this 
determinantal resultant in terms of Chern classes of $E$ and $F$. In
section \ref{LascRes}, we recall the resolution of a determinantal
variety given by A. Lascoux in \cite{Las78}. We use it in section
\ref{Compdetres}  to show that, always as for the resultant of
\cite{GKZ94}, the 
determinantal resultant can be computed as the determinant of a
certain complex. Finally in section \ref{Pdetres} we deal with the 
particular case where $X$ is a 
projective space, it is possible to be more explicit in this
situation.

This work is based on the sixth chapter of a thesis defended at the
University of Nice \cite{PhD}.


\section{Construction of the determinantal resultant}\label{secdefdetres}

Let $X$ be an irreducible projective variety over an 
algebraically closed field $\KK$.  
Let $E$ and $F$ be two vector bundles on $X$ of respective  rank 
$m$ and $n$ such that $m\geq n$.  We denote by $H=\hom(E,F)$ 
the finite dimensional vector space of   
 morphisms from $E$ to $F$. For all integer $k\geq 0$, and for all
 morphism  $\varphi \in H$, we denote by  $X_k(\varphi)$ the
 $k^{\scriptsize{th}}$-determinantal variety of $\varphi$ defined by 
$$X_k(\varphi)=\{ x\in X : \rank(\varphi(x)) \leq k \}.$$ It is well
known that this variety has 
codimension at most $(m-k)(n-k)$. Let $Y=\PP(H)$ be the projective space
associated to $H$ and let $r$ be a fixed integer such that $m \geq n >
r \geq 0$. We denote by $\nabla$ the variety 
$$\nabla=\{ \varphi \in Y \ : \ \exists x \in X \ \rank(\varphi(x)) \leq
  r \}=\{\varphi \in Y \ : \  X_r(\varphi)\neq \emptyset \},$$
and call it the \emph{resultant} variety. All along this paper, we will be
interested in this variety in the particular case where it is a
hypersurface of $Y$. 
Before stating the main theorem of this section  which gives some conditions so that
$\nabla$ is a hypersurface of $Y$, we recall briefly the
definition of the  
Hilbert scheme $\HH_P(X)$ where $P$ is a polynomial of the ring 
$\KK[\nu]$ and refer to \cite{EiHa00} for more details.

Consider the
functor $h_{P,X}$ from the opposite 
category of schemes to the category of sets,
$$h_{P,X} : \begin{array}{ccc}
            (\mathfrak{S}chemes)^{opp} & \rightarrow &
            (\mathfrak {S}ets) \\
            B & \mapsto & h_{P,X}(B),
            \end{array}$$
which associates to all scheme $B$ the set $h_{P,X}(B)$ of 
subschemes  $\XX \subset
X\times B$, flat over $B$, with fibers at each  point of $B$ 
admitting $P$ as Hilbert polynomial. This functor is 
representable and the Hilbert scheme 
$\HH_P(X)$ is defined as the scheme representing it. It
parameterizes all the subschemes of $X$ with the same Hilbert
polynomial $P$. It comes with a universal subscheme 
$\UU \subset X\times \HH_P(X)$, flat over $\HH_P(X)$ with Hilbert
polynomial $P$,  associated to the identity map. In this way, all 
subscheme $Y \subset X\times B$ flat over $B$ with Hilbert polynomial 
$P$ is isomorphic to the fiber product $Y\simeq \XX \times_{\HH_P(X)} B
\subset X\times B$ for a unique morphism $B \rightarrow
\HH_P(X)$. In what follows we will only consider the Hilbert scheme
 $\HH_2(X)$ which parameterizes all the subschemes of $X$ with 
Hilbert polynomial $P=2$, that is all the zero-dimensional subschemes
of degree 2 of $X$. 

\begin{theo}\label{defresdet} Let $r$ be a positive integer such that 
  $m \geq n > r \geq 0$, and suppose that $X$ is an irreducible projective
  variety of dimension $(m-r)(n-r)-1$. Suppose also that for all $z
  \in \HH_2(X)$ the canonical restriction morphism 
$$H \longrightarrow H^0(z,\Hom(E,F)_{|z})$$
is surjective, then $\nabla$ is an irreducible hypersurface of $Y$.
\end{theo}

\begin{preuve}  Consider the incidence variety 
  $$W=\{(x,\varphi) \in X \times Y \ : \ \rank(\varphi(x)) \leq r \}
  \subset X \times Y,$$ 
  and the two canonical projections 
  $$ X \xleftarrow{p} W \xrightarrow{q} Y.$$
  The resultant variety $\nabla$ is obtained as the projection of $W$ on
  $Y$, that is 
  $$q(W)=\{ \varphi \in Y \ : \ \exists x \in X \ \rank(\varphi(x)) \leq
  r \}=\nabla.$$ 
  Let $W_x$ be the fiber of $W$ at the point $x \in
  X$. Since, by hypothesis, $H$ generates $\Hom(E,F)$, the image of
  the evaluation  morphism at $x$ 
  $$W_x \xrightarrow{ev_x} \PP\hom(\KK^m,\KK^n)$$
  consists in all the morphisms of $\PP\hom(\KK^m,\KK^n)$ of rank
  lower or equal to $r$. These morphisms form an irreducible variety 
 of $\PP\hom(\KK^m,\KK^n)$ of codimension $(m-r)(n-r)$.  
   As the fibers of the evaluation morphism are 
  isomorphic to  
  morphisms from $E$ to $F$ which vanish at $x$, hence in particular
  irreducible, we deduce that $W_x$ is irreducible. We deduce also
  that its codimension in $H$ is the codimension of the morphisms of
  $\PP\hom(\KK^m,\KK^n)$ of rank lower or equal to $r$, that is 
  $(m-r)(n-r)$. Denoting by $s+1$ the dimension  
  of $H$, $Y$ has dimension $s$ and we obtain that  $W$
  is an irreducible projective variety of  dimension 
  $$\dim(W)=(m-r)(n-r)-1+s-(m-r)(n-r)=s-1.$$
  Now if we can show that the morphism $W 
  \rightarrow \nabla$ induced by $q$ is birational, then we show that $\nabla$ 
  is an irreducible projective variety of codimension $1$ in $H$,
  birational to $W$. Consequently, it remains  to prove that the
  morphism $W  \xrightarrow{q} \nabla$ is birational, that is its  
  generic fiber is a smooth point (notice here that since $X$ is a
  variety over a field its singular locus is of codimension at least
  one, see \cite{Hart77}, II.8.16). To do this, we will show that the
  subvariety of $Y$ consisting in all the morphism having 
  rank lower or equal to $r$ on a zero-dimensional subscheme of degree
  2 of $X$ (that is on two distinct points or on a double point) is
  of codimension greater or equal to 2.\\
  We denote by $r$ (resp. $s$) the canonical  
  projection of the universal subscheme $\UU \subset
  X\times \HH_2(X)$ of the Hilbert scheme $\HH_2(X)$ on $X$ (resp. $\HH_2(X)$),
  $$X \xleftarrow{r} \UU \xrightarrow{s} \HH_2(X).$$
  By definition of $\UU$ and $\HH_2(X)$, the morphism $s$ is flat 
and projective. As $r^*\Hom(E,F)$ is a bundle over $\UU$, it is flat
over $\UU$, and hence $r^*\Hom(E,F)$ is flat over $\HH_2(X)$.  In this
way all the 
hypothesis of the ``changing basis'' theorem are satisfied (see \cite{Hart77}
  theorem 12.11). For all $z \in \HH_2(X)$ the fiber of the morphism $s$,
  denoted $\UU_z$, is zero-dimensional in $\UU$ and hence 
  $H^1(\UU_z,r^*(\Hom(E,F))_z)=0$; we deduce that $s_*r^*\Hom(E,F)$
  is a bundle over $\HH_2(X)$ of rank $2mn$ and that its fiber at a point
  $z$ of $\HH_2(X)$ satisfies 
  $$s_*r^*\Hom(E,F)_z\simeq H^0(\UU_z,r^*\Hom(E,F)_z)\simeq
  H^0(z,\Hom(E,F)_{|z}).$$
  We consider now the incidence variety 
  $$D=\{ (z,\varphi) \in \HH_2(X) \times Y \ : \ \rank
  (s_*r^*\varphi(z)) \leq r \} \subset \HH_2(X) \times Y.$$ 
  The image of the natural projection of $D$ onto $Y$ is exactly 
  the set of $\varphi \in Y$ which vanish on a  zero-dimensional subscheme 
  of degree 2 of $X$. Let $z$ be a point of 
  the Hilbert scheme $\HH_2(X)$. It is a zero-dimensional subscheme of
  degree 2 on $X$; it is hence associated to a dimension 2 
  $\KK$-algebra, and hence of the  form
  $\Spec(\KK\oplus\KK)$ or $\Spec(\KK[\epsilon]/\epsilon^2)$.\\
  Let $D_z$ be the fiber of the projection $D\rightarrow \HH_2(X)$
  above a point $z\in \HH_2(X)$. If $z$ is of type 
  $\Spec(\KK\oplus\KK)$, that is two $\KK$-rational distinct points 
  $x$ and $y$ of $X$, then $s_*r^*\Hom(E,F)_z\simeq \Hom(E,F)_x \oplus
  \Hom(E,F)_y$. By hypothesis, the double evaluation at $x$ and $y$
  $$D_z\xrightarrow{ev_x\times ev_y} \PP\hom(\KK^m,\KK^n) \times 
  \PP\hom(\KK^m,\KK^n)$$ has image all pairs of morphisms
   $(f,g) $ such that $f(x)$ and $g(y)$ are of rank lower or equal to 
  $r$. The fibers of this morphism being isomorphic to 
  morphisms vanishing at $x$ and $y$, we deduce that the codimension
  of 
  $D_z$ in $Y$ is the one of all pairs of morphisms
  $(f,g)$ such that $f(x)$ and $g(y)$ are of rank lower or 
  equal to $r$ in $\PP\hom(\KK^m,\KK^n) \times \PP\hom(\KK^m,\KK^n)$,  
  that is $2(m-r)(n-r)$. 
  If $z$ is now of type $\Spec(\KK[\epsilon]/\epsilon^2)$,
  that is a $\KK$-rational point $x$ in $X$ and a tangent vector  
  $t \in T_xX$ at $x$, then  $s_*r^*\Hom(E,F)_z\simeq
  (\KK[\epsilon]/\epsilon^2)^{mn}$. A morphism $\varphi \in Y$
  evaluated at $x+\epsilon t$ is identified to $\varphi(x+\epsilon
  t)=f(x)+g(t)\epsilon$. As for the case of two simple points, the
  evaluation morphism 
  $$D_z\rightarrow \PP\hom(\KK^m,\KK^n) \times 
  \PP\hom(\KK^m,\KK^n) : \varphi \rightarrow (f,g)$$ has image all
  pairs of morphisms 
  $(f,g) $ such that $f(x)$ and $g(t)$ are of rank lower or equal to
  $r$ and, by the same arguments, $D_z$ is hence of 
  codimension $2(m-r)(n-r)$ in $Y$.\\
  Since the dimension of $\HH_2(X)$ is twice the dimension of $X$, we 
  deduce that 
  $$\dim(D)=2(m-r)(n-r)-2+s-2(m-r)(n-r)=s-2,$$
  and hence that the projection of $D$ onto $Y$ is of codimension
  greater or equal to $2$.
\end{preuve}

  Under the hypothesis of this theorem we define $\Res_{E,F,r}$, and
  call it the $\rth$-\emph{determinantal resultant of E and F}, to be the
  irreducible polynomial equation (defined up to a nonzero
  multiplicative constant of $\KK$) of the hypersurface $\nabla$. This
  determinantal  resultant satisfies the property~: for all $\varphi
  \in H=\hom(E,F)$, $$\Res_{E,F,r}(\varphi)=0 \Longleftrightarrow X_r(\varphi)\neq
  \emptyset.$$

\begin{rem} {\rm
In the particular case $F=\OO_X$ the $0^{\scriptsize{th}}$-determinantal
resultant  is exactly the resultant associated to a vector bundle
described in \cite{GKZ94}, chapter 3, section C. If we suppose
moreover that $E$ splits in $E=\oplus_{i=1}^m \LL_i$ where each
$\LL_i$ is a line
bundle such that its dual is very ample, the
$0^{\scriptsize{th}}$-determinantal resultant  
corresponds to the mixed resultant (\cite{GKZ94},
chapter 3, section A) associated to 
$\LL_i^*$, $i=1,\ldots,m$. Continuing the specialization, if we
suppose moreover that $X=\PP^{m-1}$ and $\LL_i\simeq \OO_X(-d_i)$,
with $d_i\geq 1$, for all $i=1,\ldots,m$, then the
$0^{\scriptsize{th}}$-determinantal resultant   
corresponds to the Macaulay resultant widely studied by
J.P. Jouanolou (see \cite{Jou91,Jou97}).
}
\end{rem}

As we can see in the hypothesis of theorem \ref{defresdet}, being
given $X$, $E$ and $F$, the $\rth$-determinantal resultant
$\Res_{E,F,r}$ not always exists. In fact such a determinantal resultant
exists only if the integer $r$ satisfies the inequality $m\geq n > r
\geq 0$ and  the equality $(m-r)(n-r)=\dim(X)+1$. From this we can see
directly that the $(n-1)^{\scriptsize{th}}$-determinantal resultant, 
which we will often call the \emph{principal determinantal
resultant}, exists if $\dim(X)=m-n$. Also we can see that if
$\dim(X)+1$ is prime  
 then the principal determinantal resultant is the only one possible
 (and exists if $\dim(X)=m-n$). 

As $\Res_{E,F,r}$ represents, when it exists, an irreducible and
reduced divisor of 
$Y=\PP(H)$, it is quite natural to ask for its degree. We begin by
fixing some notations.\\
If $V$ is a vector bundle on a variety $Z$, 
$$c_t(V)=1+c_1(V)t+c_2(V)t^2+\ldots$$
denotes its Chern polynomial. We define the Chern polynomial 
of the virtual bundle $-V$ to be 
$$c_t(-V)=\frac{1}{c_t(V)}=1-c_1(V)t+(c_1^2(V)-c_2(V))t^2+\ldots.$$
For all formal series $$s(t)=\sum_{k=-\infty}^{+\infty}c_k(s)t^k,$$ 
we set 
$$\Delta_{p,q}(s)=\det\left(
         \begin{array}{ccc}
          c_p(s) & \ldots & c_{p+q-1}(s) \\
          \vdots & & \vdots \\
          c_{p-q+1}(s) & \ldots & c_p(s)
         \end{array}
  \right).$$

\begin{prop}\label{detresdeg}
We suppose that all the assumptions of theorem \ref{defresdet} are
satisfied. Let 
$f(t,\alpha)$ be the polynomial in two variables defined by 
$$f(t,\alpha)=\sum_{k=0}^m \left(c_k(E)-(m-k+1)c_{k-1}(E)\alpha \right)t^k.$$ 
Then, the degree of the determinantal resultant $\Res_{E,F,r}$ is 
$$(-1)^{(m-r)(n-r)}\int_X\Delta_{m-r,n-r}(f(t,\alpha)/c_t(F))_\alpha ,$$
where the subscript $\alpha$ denotes the coefficient of $\alpha$  of
the  univariate polynomial $\Delta_{m-r,n-r}(f(t,\alpha)/c_t(F))$, and 
where $\int_X$ denotes the degree map of $X$.
\end{prop}

\begin{preuve} We begin again with the notations of the proof of 
  theorem \ref{defresdet}. The incidence variety 
  $$W=\{(x,\varphi) \in X \times Y \ : \ \rank(\varphi(x)) \leq r \}
  \subset X \times Y,$$ 
  where $Y=\PP(\hom(E,F))$, has the two canonical projections 
  $$ X \xleftarrow{p} W \xrightarrow{q} Y.$$
  We showed that the morphism $q$ is birational on its image which is
  exactly the variety $\nabla$ defined by the vanishing of the  
  determinantal resultant $\Res_{E,F,r}$. We have also showed that  $W
  \subset  X \times Y$ is irreducible of codimension $(m-r)(n-r)$.

  Now we denote also (abusing notations) by $p$ and $q$ both  projections
  $$ X \xleftarrow{p} X \times Y \xrightarrow{q} Y,$$
  and consider the vector bundle $p^*(\Hom(E,F))\otimes q^*(\OO_{Y}(1))$
  over $X \times Y$. The vector space of global sections of this
  vector bundle is naturally identified with $\End(H)$, the 
  vector space of endomorphisms of $H=\hom(E,F)$ (morphisms from $H$ onto itself). Let   
$$\sigma: p^*(E)\otimes q^*(\OO_{Y}(-1)) \longrightarrow 
p^*(F)$$ be the section of the vector bundle $p^*(\Hom(E,F))\otimes
q^*(\OO_{Y}(1))$ which corresponds to the identity endomorphism. Its 
$r^{\scriptsize{th}}$-determinantal variety, defined by 
$$\DD_r(\sigma)=\{ (x,f) \in X \times Y: \rank(\sigma(x,f)) \leq r
\},$$ is exactly the incidence variety $W$ since
$\sigma(x,f)=f(x)$ for all $(x,f)\in X \times Y$.
It is of codimension $(m-r)(n-r)$, which  
is its  waited codimension since the vector bundles $\EE=p^*(E)\otimes
q^*(\OO_{Y}(-1))$ and $\FF=p^*(F)$ on $X \times Y$ are 
respectively of rank $m$ and $n$. If we denote by $h$ the generic
hyperplane of $Y$ and $\alpha={q}^*(h)$, the 
Thom-Porteous formula (see \cite{ACGH85}, chapter II) and the
birationality of $q:W\rightarrow \nabla$ show that the degree of
$\Res_{E,F,r}$ is given by the coefficient of $h$ of 
$$q_*(\Delta_{n-r,m-r}(c_t(\FF-\EE)))= q_*(
(-1)^{(m-r)(n-r)}\Delta_{m-r,n-r}(c_t(\EE-\FF))),$$
that is
$$(-1)^{(m-r)(n-r)}\int_{Y} {q}_*
\left(\Delta_{m-r,n-r}(c_t(\EE)/c_t(\FF))_\alpha \right),$$
where $\int_{Y}$ denotes the degree map over $Y$.
Now  
the Chern polynomial of the bundle $\EE$ equals 
\begin{eqnarray*}
  c_t(\EE) & = & \sum_{k=0}^m \left(\sum_{i=0}^k \mathbf{C}_{m-i}^{k-i}
  c_i(p^*(E))(-\alpha)^{k-i}\right)t^k \\
& = & \sum_{k=0}^m \left(c_k(p^*(E))-(m-k+1)c_{k-1}(p^*(E))\alpha \right)t^k +
O(\alpha^2). 
\end{eqnarray*}
In this way we can see $\Delta_{m-r,n-r}(c_t(\EE-\FF))$ as a
polynomial in the variable $\alpha$. 
 As  $p^*$ commutes with Chern classes, we obtain the desired result
 from the projection formula.
\end{preuve}

\begin{rem} {\rm If we focus  on the principal case the
    preceding formula shows that the degree of a principal determinantal
    resultant is given by
    $$(-1)^{m-n+1}c_{m-n+1}\left(f(t,\alpha)/c_t(F)\right)_\alpha.$$
    If we suppose moreover that $F={\OO_X}^{\oplus n}$ then this formula
    specializes in the simple expression $(-1)^{m-n}nc_{m-n}(E)$, that
    is $nc_{m-n}(E^*)$.  
    Finally if we suppose moreover that $n=1$ (i.e. $F=\OO_X$ and
    $r=0$) we obtain $c_{m-1}(E^*)$ 
  which is well the degree of the resultant associated to
  a very ample vector
  bundle given in \cite{GKZ94}, chapter 3,
  theorem 3.10).}
\end{rem}

Classical resultants of polynomial systems are known to be
multi-homoge\-neous, depending of the geometry of the system. This
property can be stated with the notion of mixed resultant
(see \cite{GKZ94}, 
chapter 3). Let $X$ be an irreducible projective variety of dimension
$m-1$, let $\LL_1,\ldots,\LL_{m}$ be $m$ very ample line bundles
on $X$ and denote by $V_i=H^0(X,\LL_i)$. Then the mixed resultant of
$\LL_1,\ldots,\LL_{m}$ is a
polynomial over $V=\oplus_{i=1}^{m} V_i$ which is multi-homogeneous in
each $V_i$ of degree $\int_X \prod_{j\neq i} c_1(\LL_j)$ (recall that
the mixed resultant can be seen as the principal determinantal
resultant of $E=\oplus_{i=1}^{m}\LL_i^*$ and $F=\OO_X$). In the
context of the determinantal resultant this multi-homogeneity 
appears as follows.

\begin{prop}\label{multideg}
We suppose that all the assumptions of theorem \ref{defresdet} are
satisfied and that $E=\oplus_{i=1}^mE_i$ where each $E_i$ is a line
bundle on $X$. Let 
$f(t,\alpha_1,\ldots,\alpha_m)$ be the polynomial in $m+1$ variables
defined by  
$$f(t,\alpha_1,\ldots,\alpha_m)=\prod_{i=1}^m (1+(c_1(E_i)-\alpha_i)t).$$
The $\rth$-determinantal resultant is multi-homogeneous in each vector
space 
$H_i=\hom(E_i,F)$, $i=1,\ldots,m$, and its degree with respect to $H_i$ is 
$$(-1)^{(m-r)(n-r)}\int_X\Delta_{m-r,n-r}(f(t,\alpha_1,\ldots,\alpha_n)
/c_t(F))_{\alpha_i} ,$$ 
where the subscript $\alpha_i$ denotes the coefficient of $\alpha_i$
of the multivariate polynomial $\Delta_{m-r,n-r}(f(t,\alpha_1,\ldots,\alpha_n)
/c_t(F))$ (in variables $\alpha_1,\ldots,\alpha_m$). 
\end{prop}

\begin{proof}
This proof is just a refinement of the proof of proposition
\ref{detresdeg}. If the vector bundle $E$ splits into
$E=\oplus_{i=1}^mE_i$ then clearly $H=\hom(E,F)=\oplus_{i=1}^m
\hom(E_i,F)=\oplus_{i=1}^mH_i$. Each morphism $\varphi \in H$ can be
decomposed as $\varphi=\oplus_{i=1}^m \varphi_i$ where each $\varphi_i
\in H_i$. Now multiplying each morphism $\varphi_i$ by its own nonzero constant
$\lambda_i$ do not change the rank, that is
$$\forall x \in X \ \
\rank(\varphi(x))=\rank(\oplus_{i=1}^m\lambda_i\varphi_i(x)).$$
It follows that $\Res_{E,F,r}$ is multi-homogeneous with respect to
each $H_i$. Consequently we can see our determinantal resultant in
$Y_1\times\ldots\times Y_m$, where $Y_i$ denotes the projective space
$\PP(H_i)$, instead of $Y=\PP(H)$.

Taking again the proof of proposition \ref{detresdeg} the incidence
variety $W$ is obtained as the zero locus of the canonical section 
$$\sigma : \EE=\bigoplus_{i=1}^m \left( p^*(E_i)\otimes q^*(\OO_{Y_i}(-1)) \right) \rightarrow \FF=p^*(F).$$
Denote by $h_i$ the generic hyperplane of $Y_i$ for all
$i=1,\ldots,m$. The Chern polynomial of the bundle $\EE$ can be written
$$c_t(\EE)=\prod_{i=1}^{m}(1+(c_1(E_i)-\alpha_i)t),$$
where $\alpha_i=q^*(h_i)$, and we obtain the desired result as in
proposition \ref{detresdeg}.
\end{proof}


\section{Resolution of a determinantal variety}\label{LascRes}

In this section we will focus on the following problem: being given 
two vector bundles $E$ and $F$ on a scheme $X$ and a sufficiently
generic morphism $\varphi : E \rightarrow F$,
describe the locus $Y$ of points of $X$ where $\varphi$ has rank lower or
equal to a given integer $r$; more precisely, define an exact complex 
 $L(\varphi,r)^\bullet$ of vector bundles over $X$ which gives a
 resolution of the trivial vector bundle of $Y$, $\OO_{Y}$. This 
 problem has been solved by Alain Lascoux in its paper titled {\it ``Syzygies 
  des vari\'et\'es d\'eterminantales``}, \cite{Las78}. Its result is
based on the fact that  $Y$ is birational to a subscheme  
 $Z$ of a grassmannian $G$ on $X$,
that is there exists a commutative diagram
$$\begin{CD}    
Z @>>> G \\ 
@VVV  @V\pi VV \\ 
Y @>>> X,\\  
\end{CD}$$
such that the restriction of the map $\pi$ to $Z$ is birational on $Y$, and
such that $\OO_Z$ admits a known resolution $K^\bullet$ on $G$ (this
resolution is the Koszul complex of a section of a vector bundle over
$G$). The author then studies the spectral sequence associated to the
projection of the complex $K^\bullet$ on $X$, and constructs in
this way a resolution of $\OO_Y$ on $X$. Notice that we have to
suppose here that the algebraically closed field $\KK$ is
of characteristic zero, hypothesis needed in
\cite{Las78}. However we 
will see at the end of this section that this hypothesis is not
needed in the particular case of interest
$r=\min(\rank(E),\rank(F))-1$ (the principal case).

In what follows we gather (from the original paper of Alain Lascoux
\cite{Las78}) the principal ingredients in order to state 
the nice results obtained by Alain Lascoux we will use in the next
section to compute the determinantal resultant. Before giving the 
resolution $L(\varphi,r)^\bullet$ itself 
 we recall first some definitions on  partitions and 
Schur functors.

\subsection{Partitions and Schur functors}

\subsubsection{Partitions} Let $r$ be an integer. A \textit{partition}
is a $r$-uples $I=(i_1,\ldots,i_r)$ in $\NN^r$ such that $0\leq i_1\leq
\ldots \leq i_r$; its \emph{length} is the number of its nonzero elements,
and its \emph{weight}, denoted $|I|$, is $i_1+\ldots +i_r$. We will identify
naturally two partitions which differ only by adding zeros on the
left, but however it will be useful to consider the partition 
$(0,\ldots,0)\in \NN^r$ which will be denoted $O_r$. 

A partition is usually represented by an interval of 
$\NN\times \NN$, its \textit{Ferrer diagram}; the first line  
contains $i_r$ boxes, the second  $i_{r-1}$ boxes, and so on...   
This yields the following picture~:
\begin{center}
\includegraphics[width=3cm,height=3cm]{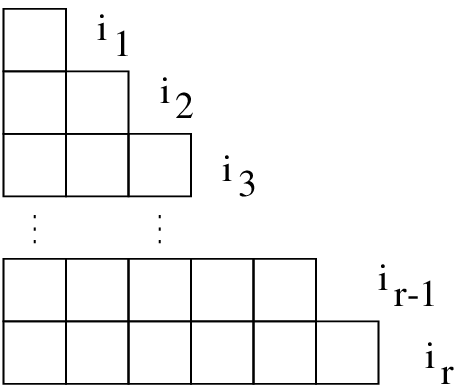}
\end{center}
A partition $J$ is said to be \textit{contained} in $I$, denoted
$J\subseteq I$, if it is true for their Ferrer diagrams.

On the set of all partitions we define an involution map which sends
each partition $I$ to the partition
$I^\star$, called the \emph{dual} partition of $I$, such that its
Ferrer diagram is the transposed Ferrer diagram of $I$. For
instance $(1,2,4)^\star=(1,1,2,3)$.

If $I$ and $J$ are two partitions $(i_1,\ldots,i_r)$ and 
$(j_1,\ldots,j_q)$ such that $i_r\leq j_1$, the \emph{concatenation}
$IJ=conc(I,J)=(i_1,\ldots,i_r,j_1,\ldots,j_q)$ is a partition. It is
possible to 
extend this operator as follows. Let $r$ and $q$ be two integers, $I\in
\ZZ^r$, $J\in \ZZ^q$, and consider the set 
$\{i_1,i_2+1,\ldots,i_r+r-1,j_1+r,\ldots,j_q+r+q-1\}$. If all these
numbers are positive and distinct, this set can be write in a unique
way as $\{h_1,\ldots,h_{r+q}+r+q-1 \}$, with 
$H=(h_1,\ldots,h_{r+q})$ a partition. We say that $H$ is the 
\textit{concatenation} of $I$ and $J$, with \textit{ampleness} $n(I,J)$,
the minimal number of transpositions to rectify the first set onto the
second. Otherwise we set $conc(I,J)=\emptyset$ and 
$n(I,J)=\infty$. 

\subsubsection{Tablo\"{\i}ds - Young diagrams} 
Let $E$ be a set and $I$ be a partition. A \textit{tablo\"{\i}d} of 
diagram $I$ with values in $E$ is a filling of the diagram of $I$ with
elements of $E$.

Let $\tau$ be a tablo\"{\i}d and $\mu$ be a permutation of the boxes
of the diagram of $I$. We denote by $\tau^\mu$ the tablo\"{\i}d
induced by this permutation and we define $Tabl_IE$ as the free $\ZZ$-module
with basis all tablo\"{\i}ds of diagram $I$.

Suppose now  the set $E$ has a total ordering.  
A \textit{Young diagram} of a diagram $I$ is a tablo\"{\i}d of 
diagram $I$ such that each column read from the low to the top 
is a strictly growing sequence, and each line, from
the left to the right, is a growing sequence.

\subsubsection{Schur functors} 
Let $E$ be a vector space, $r$ be an integer and $I$ be an element of 
$\NN^r$. We denote by $S^IE$ the tensor product $S^{i_1}E
\otimes \ldots \otimes S^{i_r}E$, where $S^\bullet E$ is the symmetric
algebra of $E$. In the same way $\wedge^IE=\wedge^{i_1}E\otimes \ldots
\otimes \wedge^{i_r}E$, where $\wedge^\bullet E$ is the exterior
algebra of $E$. 

It is possible to represent the elements of $S^iE$ by line
tablo\"{\i}ds, 
and the elements of $\wedge^iE$ by column tablo\"{\i}ds (such a
tablo\"{\i}d is not unique in general). Let $I$ be a 
partition and $J=I^\star$ its dual, we have canonical morphisms 
$$ \psi_S:Tabl_IE \rightarrow S^IE$$
$$ \psi_\wedge:Tabl_IE \rightarrow \wedge^JE$$
which are obtained by reading the lines or the columns of a
tablo\"{\i}d of diagram $I$.

Let $\rho_I^\prime$ be the endomorphism of $Tabl_IE$
$$ \begin{array}{ccc}
              Tabl_IE & \rightarrow & Tabl_IE \\
              \tau & \mapsto & \sum_{\mu} (-1)^\mu \tau^\mu \\
              \end{array}$$
where the sum is running over all permutations $\mu$ preserving 
the columns of the diagram, with
$(-1)^\mu=\textrm{signature}(\mu)$. Then there exists a unique
morphism which makes the following diagram commutative~:
$$\begin{CD}    
Tabl_IE @>\rho_I^\prime>> Tabl_IE \\ 
@V\psi_{\wedge}VV  @V\psi_{S}VV \\ 
\wedge^JE @>\rho_{I}>> S^IE.\\  
\end{CD}$$
Choosing a total ordering on a basis of $E$, one shows that the image
of the map $\rho_{I}$ is a vector space with basis the images 
of Young diagrams by the composed map $\psi_{S} \circ \rho_I^\prime$. 
 
This construction for a vector space 
$E$ is easily extended to a vector bundle  
$E$ on a scheme $X$. We define the \textit{Schur functor of
  index $I$} (from the category of vector bundles on $X$ to itself),
denoted  indifferently $S_I$ or 
$\wedge_J$, as the image functor  
of $\wedge^J$ in $S^I$ by $\rho_I$ (we set 
$S_\emptyset=0$).\\ 
Notice that if the length of $I$ is strictly greater 
than the rank of $E$ then $S_IE=0$, since 
$\wedge^{length(l)}E=0$.

\subsection{The  resolution}\label{seclascres}

Let $E$ and $F$ be two vector bundles of respective rank 
$m$ and $n$, $m\geq n$, on a scheme 
$X$ and $\varphi$ a morphism between them $\varphi : E \rightarrow
F$. Let also $r$ be a positive integer strictly lower to $n$ and 
$Y$ be the subscheme of $X$ defined by 
$$Y=\{ x \in X : \rank(\varphi(x)) \leq r\}.$$
Its trivial bundle $\OO_Y$ is defined as 
$\coker(\varphi')$ where 
$$ \varphi' : \wedge^{r+1}E \otimes \wedge^{r+1} F^*
\rightarrow \OO_X$$
is associated to the morphism $\wedge^{r+1}\varphi$.\\
Consider now the relative grassmannian 
$G=G_r(F)\xrightarrow{\pi} X$ of quotients of $F$ of rank $q=n-r$, and
$$0 \rightarrow S \rightarrow \pi^*(F) \rightarrow Q \rightarrow 0$$
its tautological exact sequence. Let $Z$ be the subscheme of $G_r(F)$
defined as the zero locus of the composed morphism $\pi^*E \rightarrow \pi^*F
\rightarrow Q$, i.e. $\OO_Z$ is defined as the cokernel 
of the induced map $\psi:\pi^*E \otimes Q^* \rightarrow \OO_{G}$. We
suppose that $Z$ is {\it locally a complete intersection}, that is 
 the Koszul complex associated to $\psi$,
$$0\rightarrow K^{-m(n-r)}\rightarrow \ldots \rightarrow
K^{p}=\wedge^{-p}(\pi^*E \otimes Q^*) \rightarrow \ldots \rightarrow
K^0=\OO_G,$$ 
is exact (this hypothesis implies that the restriction of $\pi$ to 
$Z\rightarrow Y$ is birational, see \cite{Las78}). We can then construct a complex 
 $L(\varphi,r)^\bullet$ from the non-degenerated spectral sequence 
$$ \EE^{s,t}=R^t\pi_*(\wedge^{-s}(\pi^*E \otimes Q^*)),$$
by setting $L^p=\bigoplus_{s+t=p}\EE^{s,t}$.\\
Denoting $n(I)=n(I,O_r)$ and $I^\prime=conc(I,O_r)$ the Cauchy formula
and the Bott theorem (see \cite{Las78}) show that 
$$L^p=\bigoplus_{s+t=p}\EE^{s,t}=\bigoplus_{-|I|+n(I)=p}\wedge_IE
\otimes S_{I^\prime} F^*,$$  
where $S_{I^\prime} F^* = R^{n(I)}\pi_*(S_IQ^*)$ (which 
implies that $I$ must be of length lower or equal to $q$ so that  
$S_{I^\prime} F^*$ is non-zero). The following lemma gives an explicit
description of $I^\prime$ and $n(I)$~: 
\begin{lem}\label{comblascoux} {\rm (\cite{Las78}, 5.10)} Let
  $I=(i_1,\ldots,i_q)$ and let 
  $p(I)$ be the dimension of the greatest square contained in the
  diagram of $I$, then~:
  \begin{itemize}
    \item $p=p(I)$ is such that $i_{q-p}+q-p-1\leq n-1$ and 
      $i_{q-p+1}+q-p\geq n$,
     \item If $i_{q-p+1} < p+r$ then $I^\prime=\emptyset$,
     \item If $i_{q-p+1} \geq p+r$ then $n(I)=pr$ and 
       $$I^\prime=(i_1,\ldots,i_{q-p},\underbrace{p,\ldots,p}_{r},i_{q-p+1}-r,\ldots, i_q-r).$$

   \end{itemize}
\end{lem}

Now we have to link the terms $L^p$ to obtain the complex 
$L(\varphi,r)^\bullet$. The lemma 5.13 of \cite{Las78} tells us that
if $(I,H)$ is a pair of 
partitions such that $H\subset I$, $H^\prime\neq \emptyset$,
$I^\prime\neq \emptyset$ and $|I|-n(I)=|H|-n(H)-1$, then $I$ and $H$
have the same respective parts, excepted one, and the same is true for
${I^\prime}^\star$ and ${H^\prime}^\star$. More precisely we have 
$$I=\bullet - \bullet \, i \bullet - \bullet, \
{I^\prime}^\star=\circ - \circ \, i' \circ - \circ,$$
$$H=\bullet - \bullet \, h \bullet - \bullet, \
{H^\prime}^\star=\circ - \circ \, h' \circ - \circ,$$
where $\bullet - \bullet$ (resp. $\circ - \circ$) is the set where 
$I$ and $H$ are the same (resp. ${I^\prime}^\star$ and
${H^\prime}^\star$). For such a pair ($I,H$) we can hence define the 
morphism
$$\psi_{I,H}:\wedge_IE \otimes S_{I^\prime}F^* \rightarrow \wedge_HE
\otimes S_{H^\prime}F^*$$
as the composition of the canonical injective and surjective maps 
(coming from the Pieri formula, see \cite{Las78} (1.5.3)) and the
contraction morphism 
$$d_s:\wedge^iE\otimes \wedge^{i'}F\rightarrow \wedge^{i-h}E\otimes
\wedge^hE\otimes \wedge^{i'-h'}F^*\otimes\wedge^{h'}F^* \rightarrow 
\wedge^hE\otimes \wedge^{h'}F,$$ (with $s=i-h=i'-h'$)~:
$$\wedge_IE\otimes \wedge_{I^\prime}F\longrightarrow 
\wedge_{\bullet - \bullet}E\otimes \wedge^iE\otimes \wedge^{i'}F^* 
\otimes \wedge_{\circ - \circ}F^*$$
$$\xrightarrow{d_s} \wedge_{\bullet - \bullet}E\otimes \wedge^hE\otimes
\wedge^{h'}F^* \otimes \wedge_{\circ - \circ}F^* \longrightarrow
\wedge_HE\otimes \wedge_{H^\prime}F.$$ 
For all pair $(I,H)$ such that $H$ is not contained in $I$ we set 
$\psi_{I,H}=0$. In this way we define the differential maps of the
complex $L(\varphi,r)^\bullet$ by~:
$$L^p=\bigoplus_{-|I|+n(I)=p}\wedge_IE
\otimes S_{I^\prime} F^* \xrightarrow{\psi_p=\oplus\psi_{I,H}}
L^{p+1}=\bigoplus_{-|H|+n(H)=p+1}\wedge_HE \otimes S_{H^\prime} F^*.$$
The complex $L(\varphi,r)^\bullet$ that we obtain is such that 
$L^0=\OO_X$, and such that the term on the far left corresponds to the
partition $I=(\underbrace{m,\ldots,m}_{q})$ for which $n(I)=qr$,
that is $L^{qr-mq}$. We can now state the following theorem~:

\begin{theo}\label{lascres} {\rm \cite{Las78}} Suppose that the
  subscheme $Z$ of  
  $G_r(F)$ is locally a complete intersection, then the complex
  $L(\varphi,r)^\bullet$ is a minimal resolution of $\OO_Y$.
\end{theo}

Before ending this section we would like to say a little more on the morphism
$L^{-1}\rightarrow \OO_X$ and on the principal case.

\noindent $\bullet$ The morphism $L^{-1}\rightarrow \OO_X$.\\
Lemma \ref{comblascoux} describes all  pairs of partitions
$(I,I^\prime)$ such that $|I|-n(I)=1$. In fact if we want  $I^\prime$ 
to be zero we need that $i_{p-q+1}\geq p+r$ and hence we obtain $|I|\geq
p(p+r)$. As in this case $n(I)=pr$, the unique partition $I$ such that
$|I|-n(I)=1$ and $I^\prime\neq 0$ is the line partition $I=(r+1)$. Its
dual partition $I^\prime$ is then the column partition 
$I^\prime=(\underbrace{1,\ldots,1}_{r+1})$, and hence we deduce that 
$L^{-1}=\wedge^{r+1}E\otimes\wedge^{r+1}F^*$. The morphism 
$L^{-1} \rightarrow L^0=\OO_X$ is then reduced, by the
description of the complex, to the contraction morphism 
$\wedge^{r+1}E\otimes\wedge^{r+1}F^* \rightarrow \OO_X$ associated to
the morphism $\wedge^{r+1}\varphi$. We can hence check that the
cokernel of $L^{-1}\rightarrow \OO_X$ is well $\OO_Y$ as desired.

\noindent $\bullet$ The principal case.\\
We focus here on the case where $r=n-1$, i.e. $q=1$. Let $s$ be a
fixed integer such that $1\leq s \leq m-n+1$. We want to explicit all  
pairs of partitions $(I,I^\prime)$ such that $I\neq\emptyset$,
$I^\prime \neq \emptyset$ and $|I|-n(I)=s$. As $q=1$, $I$ is 
always a line partition $(i_1)$ and hence $p(I)=1$. Moreover, if we
want $I^\prime$ to be zero, we must have  $i_1\geq n$ and then 
$n(I)=n-1$. We deduce that $I=(n+s-1)$ and  
$I^\prime=(\underbrace{1,\ldots,1}_{r-1},s)$. We obtain in this way 
$$L^{-s}=\wedge^{n+s-1}E\otimes S^{s-1}F^*\otimes \wedge^{n}F^*.$$
We can then easily see that the complex we obtain in this case 
 is the well known Eagon-Northcott complex. Consequently in this case
 we do not have to suppose 
 that the ground field $\KK$ is of characteristic zero since the
 Eagon-Northcott complex gives a resolution of $\OO_Y$ even if $\KK$
 is not of characteristic zero (see for instance \cite{BrVe80}).


\section{Computation of the determinantal resultant}\label{Compdetres} 
In section \ref{secdefdetres} we have defined the 
$\rth$-determinantal resultant which gives a necessary and sufficient
condition so that a given morphism between two vector bundles on an
irreducible projective variety  is of rank lower or equal to the
integer $r$ in at least one point. We will now show that it is
possible to give an explicit representation of the determinantal
resultant $\Res_{E,F,r}$ (when it exists), more precisely we will show
that this resultant is the determinant of a certain complex. The key
point here is to obtain a resolution of the incidence variety,
resolution that we will get from the preceding section, and to project
it on the space of parameters, that is on $Y=\PP(\hom(E,F))$.

$X$ is always an irreducible and projective variety of dimension 
$(m-r)(n-r)-1$ over an algebraically closed field $\KK$ of
characteristic zero, where $m,n$ and 
$r$ are three positive integers such that 
$m\geq n > r \geq 0$. Let $E$ and $F$ be two vector bundles on $X$ of
respective rank $m$ and $n$. We suppose that for all $z \in \HH_2(X)$
the restriction morphism 
$H \longrightarrow H^0(z,\Hom(E,F)_{|z})$ is surjective. Then, 
by theorem \ref{defresdet} the $\rth$-determinantal 
resultant $\Res_{E,F,r}$ exists; for all 
$\varphi \in Y=\PP(H)$, it satisfies 
$$\Res_{X,E,F}(\varphi)=0 \Longleftrightarrow X_r(\varphi)\neq
\emptyset,$$
and define a divisor  $\nabla$ on $Y$.

We denote by $p$ and $q$ both projections
$$ X \xleftarrow{p} X \times Y \xrightarrow{q} Y,$$
and we consider the incidence variety 
$$W=\{(x,\varphi) \in X \times Y \ : \
\rank(\varphi(x)) \leq r \} \subset X \times Y.$$ 
In the proof of proposition \ref{detresdeg}, we have seen that the
canonical section 
$$\sigma: p^*(E)\otimes q^*(\OO_{Y}(-1)) \longrightarrow 
p^*(F)$$
of the vector bundle $p^*(\Hom(E,F))\otimes q^*(\OO_{Y}(1))$
on $X \times Y$ is such that its
$\rth$-determinantal variety 
defined by  
$$\DD_r(\sigma)=\{ (x,\varphi) \in X \times Y: \rank(\sigma(x,\varphi)) \leq r
\},$$ is exactly the incidence variety $W$, which is of 
codimension $(m-r)(n-r)$. This section can be seen as the 
``universal section'' since the restriction of  $\sigma$ to 
 a fiber $X\times\{\varphi \}$ is just the section $\varphi \in Y$.
 To simplify the notations we set $E'=p^*(E)\otimes q^*(\OO_{Y}(-1))$
 and $F'=p^*(F)$, which are two vector bundles over $X \times Y$ 
of respective rank $m$ and $n$.

We consider now the relative grassmannian 
$$\pi:G_r(F')\rightarrow X \times Y,$$ 
of quotients of rank $n-r$, which is isomorphic to $G_r(F)\times
Y$. Denoting 
$$0\rightarrow S\rightarrow \pi^*(F')\rightarrow Q \rightarrow 0$$ 
its tautological exact sequence, we define the closed subscheme 
$Z$ of $G_r(F')$ as the zero locus of the composed morphism ~:
$$\sigma^\sharp \ : \ \pi^*(E') \xrightarrow{\pi^*(\sigma)} \pi^*(F') \rightarrow Q.$$
The morphism $\sigma^\sharp$ is hence a section of the vector bundle 
$\pi^*(E')\otimes Q^*$ of rank $m(n-r)$. We obtain the following
commutative diagram
$$\begin{CD}    
Z @>>>  G_r(F)\times Y \\
  @V \pi_{|Z} VV          @V\pi VV \\ 
\DD_r(\sigma) @>>> X \times Y,\\  
\end{CD}$$
where the two horizontal arrows are the canonical injections.
 The morphism $\pi_{|Z}$, restriction of $\pi$ to $Z$ on $\DD_r(\sigma)$, is 
birational and hence  $Z$ is of codimension $m(n-r)$, that is 
of codimension the rank of the vector bundle $\pi^*(E')\otimes
Q^*$. If we show that 
$Z$ is locally a complete intersection, i.e.  
that the Koszul complex associated to $\sigma^\sharp$,
$$0\rightarrow K^{-m(n-r)}\rightarrow \ldots \rightarrow
K^{p}=\wedge^{-p}(\pi^*E \otimes Q^*) \rightarrow \ldots \rightarrow
K^0=\OO_{G_r(F')},$$
is exact, then the theorem \ref{lascres} gives us a 
resolution of $\OO_W$ by vector bundles over $X\times
Y$, which is the complex $L(\sigma,r)^\bullet$ (see section 
\ref{seclascres}).

\begin{lem} Under the  hypothesis of theorem \ref{defresdet}, 
$Z$ is locally a complete intersection.
\end{lem}
\begin{preuve}
As the variety $Y$ is smooth, we deduce that the projection 
$\tau: G_r(F)\times Y \rightarrow G_r(F)$ is a smooth 
morphism. The restriction $\tau_{|Z}:Z\rightarrow G_r(F)$ is hence 
a flat morphism. Its geometric fiber at the point $(x,V_x)$,
where $V_x$ is a vector subspace of dimension $r$ of $F_x$,
is identified to the set of morphisms $f\in Y$ such that 
$\Im(f(x))\subset V_x$; this fiber is hence smooth. It follows
that the morphism $\tau_{|Z}:Z\rightarrow G_r(F)$ is smooth (see 
\cite{Hart77}, theorem III,10.2). Moreover the restriction of 
$\sigma^\sharp$ at each of its geometric fibers is transversal to the
zero section; as the codimension of 
$Z$ equals the rank of $\pi^*(E')\otimes Q^*$ we deduce that 
the Koszul complex associated to $\sigma^\sharp$ is exact (see 
\cite{GKZ94}, proposition II,1.4), and hence that $Z$ is locally 
a complete intersection.
\end{preuve}

The complex $L(\sigma,r)^\bullet$ is hence a 
resolution of $\OO_W$ over $X\times Y$. It is of the form
$$\ldots \rightarrow \bigoplus_{-|I|+n(I)=p}\wedge_IE'\otimes
 S_{I^\prime} {F'}^* \rightarrow \ldots \rightarrow
 \OO_{X\times Y}.$$

Let $\MM$ be a line bundle over $X$. We denote by
$\mathcal{L}^\bullet$ the 
complex $L(\sigma,r)^\bullet\otimes p^*(\MM)$; it is a 
resolution of $\OO_W\otimes p^*(\MM)$. The complex 
$\mathcal{L}^\bullet$ of $\OO_{X\times 
   Y}$-modules induces the complex ${q}_*(\mathcal{L}^\bullet)$
of $\OO_{ Y}$-modules. Its term on the far right  is 
$${q}_*(\mathcal{L}^0)={q}_*(\OO_{X\times Y}\otimes
p^*(\MM))=\OO_{ Y}\otimes H^0(X,\MM),$$
and for all $p=-1,\ldots,-(m-r)(n-r)$, we have 
$${q}_*(\mathcal{L}^p) =  {q}_*\left(\bigoplus_{-|I|+n(I)=p}
  \wedge_I\left(p^*(E)\otimes q^*(\OO(-1))\right) \otimes
  S_{I'}(p^*(F))\otimes p^*(\MM)\right)$$
$$ =  \bigoplus_{-|I|+n(I)=p} {q}_*\left(\wedge_I(p^*(E))\otimes
  {q^*(\OO(-1))}^{\otimes |I|} \otimes S_{I'}(p^*(F))\otimes p^*(\MM)\right) $$
$$ =  \bigoplus_{-|I|+n(I)=p} {q}_*\left( p^*(\wedge_I(E)\otimes
  S_{I'}(F)\otimes \MM ) \otimes {q^*(\OO_{ Y}(-1))}^{\otimes |I|}\right)$$
$$ =  \bigoplus_{-|I|+n(I)=p} H^0(X,\wedge_I(E)\otimes
  S_{I'}(F)\otimes \MM)\otimes \OO_{ Y}(-|I|).$$
 
\begin{rem} { \rm The Schur functors definition we gave show
    immediately that if $I$ is a 
    partition, $E$ a vector bundle and $L$ a line bundle, 
    then $\wedge_I(E\otimes L) \simeq L^{\otimes |I|}\otimes \wedge_I(E)$.}
\end{rem}

\noindent We come to the following definition~:

  \begin{defn} We will say that the line bundle $\MM$
    \textit{stabilizes} the  
complex ${q}_*(\mathcal{L}^\bullet)$ if for all 
$p=0,\ldots,-(m-r)(n-r)$, and all partition $I$ such that 
$-|I|+n(I)=p$, all the cohomology groups $H^i(X,\wedge_I(E)\otimes
  S_{I'}(F)\otimes \MM)$ vanish for all $i>0$.
\end{defn}

  We will now show that if $\MM$ stabilizes the complex 
  ${q}_*(\mathcal{L}^\bullet)$ then this complex  
  is acyclic. To do this consider an injective  
  resolution $I^{\bullet \bullet}$ of the complex
  $\mathcal{L}^\bullet$. We can choose this resolution

  $$\begin{array}{ccccccccccc}
    & & 0 & & & & 0 & & 0 & &   \\
    & & \uparrow & & & & \uparrow & & \uparrow & &   \\
    0 & \rightarrow & I^{-(m-r)(n-r),s}  
  & \rightarrow & \ldots & \rightarrow & I^{0,s} & \rightarrow & J^s & \rightarrow & 0 \\
    & & \uparrow & & & & \uparrow & & \uparrow & &   \\
    & & \vdots & & & & \vdots & & \vdots & &   \\
    & & \uparrow & & & & \uparrow & & \uparrow & &   \\
    0 & \rightarrow & I^{-(m-r)(n-r),0}  
  & \rightarrow & \ldots & \rightarrow & I^{0,0} & \rightarrow & J^0 & \rightarrow & 0 \\
    & & \uparrow & & & & \uparrow & & \uparrow & &   \\
   0 & \rightarrow & \mathcal{L}^{-(m-r)(n-r)}  
  & \rightarrow & \ldots & \rightarrow & \mathcal{L}^0 & \rightarrow &
  \OO_W\otimes p^*(\MM) & \rightarrow  & 0 \\
   & & \uparrow & & & & \uparrow & & \uparrow & &   \\
   & & 0 & & & & 0 & & 0 & &   \\
   \end{array}$$

such that~:

\begin{itemize}
\item the complex $J^\bullet$ is an injective resolution of 
  $\OO_W\otimes p^*(\MM)$,
\item each column is an exact complex,
\item each line is also an exact complex (since 
  $\mathcal{L}^\bullet$ is a resolution of  $\OO_W\otimes
  p^*(\MM)$).
\end{itemize}

There exists two spectral sequences which converge to the 
hyper-direct image of the complex $\mathcal{L}^\bullet$,
$\RR^i{q}_*(\mathcal{L}^\bullet)$, which is by definition the 
$i^{\scriptsize{th}}$ cohomology sheaf of the total complex 
associated to the double complex ${q}_*(I^{\bullet \bullet})$. These 
two spectral sequences 
correspond respectively to column and line filtrations of this double
complex ${q}_*(I^{\bullet \bullet})$. 

The first spectral sequence is of the form 
$$ E_1^{p,q}=R^q{q}_*(\mathcal{L}^p) \Rightarrow
\RR^i{q}_*(\mathcal{L}^\bullet).$$
The cohomology of the columns of the double complex 
$q_*(I^{\bullet\bullet})$ vanishes for $q>0$ since we have 
supposed that $\MM$ stabilizes the complex 
${q}_*(\mathcal{L}^\bullet)$, which implies that 
$R^i{q}_*(\mathcal{L}^j)=0$ for $i>0$ (using  proposition
4.2.2 of \cite{Pot97}). This spectral sequence  hence 
degenerates at the second step and we obtain 
$$ E_2^{p,q}=0 \ {\rm if} \  q \neq 0, \ \
E_2^{p,0}=H^p({q}_*(\mathcal{L}^\bullet)),$$ 
where $E_2^{p,0}$ vanishes for $p>0$. It follows 
$\RR^i{q}_*(\mathcal{L}^\bullet)=0$ if $i>0$, and 
$\RR^i{q}_*(\mathcal{L}^\bullet)=H^i({q}_*(\mathcal{L}^\bullet))$
if $i\leq 0$. \\
The second spectral sequence corresponds to the filtration of the double
complex $q_*(I^{\bullet \bullet})$ by lines. Each complexes 
$$ 0  \rightarrow  I^{-(m-r)(n-r),i}  
   \rightarrow  \ldots  \rightarrow  I^{0,i}  \rightarrow  J^i  \rightarrow  0 $$
is an exact complex of injective sheaves. We deduce that this complex
splits (if we have an injective morphism $0\rightarrow I
\xrightarrow{s} K$, with 
$I$ an injective object, then the identity $Id:I \rightarrow I$ can be
lifted to a morphism $p:K\rightarrow I$ such that $p\circ s=Id$). 
Applying the functor ${q}_*$ to this complex gives hence an exact split
complex. Our spectral sequence, at the first step, is hence such that 
$${E'}_1^{0,q}={q}_*(J^q), \ {\rm if} \ {E'}_1^{p,q}=0 \ {\rm
  if } \ p\neq 0.$$ 
It degenerates at the second step, we have 
$$ {E'}_2^{p,q}=0 \ {\rm if} \  p \neq 0, \ {\rm and} \ 
{E'}_2^{0,q}=R^q{q}_*(\OO_W\otimes  p^*(\MM)).$$
We deduce $\RR^i{q}_*(\mathcal{L}^\bullet)=0$ if $i<0$, 
and $\RR^i{q}_*(\mathcal{L}^\bullet)=R^i{q}_*(\OO_W\otimes  p^*(\MM))$
if $i\geq 0$. 

The comparison of our two spectral sequences show that 
$H^i({q}_*(\mathcal{L}^\bullet))=0$ if $i<0$, and that  
$$H^0({q}_*(\mathcal{L}^\bullet))={q}_*(\OO_W\otimes  p^*(\MM))=H^0(X,\MM)\otimes {q}_*(\OO_W).$$

The complex ${q}_*(\mathcal{L}^\bullet)$ is hence acyclic and its
cohomology in degree 0 is  $H^0(X,\MM)\otimes
{q}_*(\OO_W)$. We consider now its associated graded complex which is
a complex of 
$S^\bullet(H^*)$-modules, that we denote 
$\mathcal{C}(\mathcal{\MM})^\bullet$, and which is such that 
$$\mathcal{C}(\mathcal{\MM})^0=H^0(X,\MM)\otimes S^\bullet(H^*), \ {\rm and},$$
$$C(\MM)^{p}= \bigoplus_{-|I|+n(I)=p} H^0(X,\wedge_I(E)\otimes 
  S_{I'}(F)\otimes \MM)\otimes S^\bullet(H^*).$$
We denote also by $\mathbf{C}(H)$ the fraction field of 
 $S^\bullet(H^*)$. We have the following theorem~:
\begin{theo}\label{detresdet} Suppose the hypothesis of theorem
  \ref{defresdet} are satisfied. If $\MM$ stabilizes the complex
  $q_*(\LL^\bullet)$ then  
$$\det(\mathcal{C}(\mathcal{\MM})^\bullet \otimes \mathbf{C}(H)
)=\Res_{E,F,r}.$$
\end{theo}
\begin{preuve}
The morphism $q:W\rightarrow \nabla$ being birational,    
the complex $\mathcal{C}(\mathcal{\MM})^\bullet$ is 
generically exact, that is  
the complex of vector spaces (of finite dimensions) 
$\mathcal{C}(\mathcal{\MM})^\bullet \otimes \mathbf{C}(H)$ is exact.
 By the theorem 30 of \cite{GKZ94} we deduce that  
$$\det(\mathcal{C}(\mathcal{\MM})^\bullet \otimes \mathbf{C}(H)
)=\sum_{D} \left(\sum_{i} 
(-1)^i\mathrm{mult}_{(Z_D)}(H^i(\mathcal{C}(\mathcal{\MM})^\bullet)\right).D,$$
where the sum is running over all the irreducible polynomials $D$ of 
$S^\bullet(H^*)$, and where  
$\mathrm{mult}_{(Z_D)}$ denotes the multiplicity along the irreducible
hypersurface associated to the polynomial $D$. 
The results of \cite{Ser55} show that all the modules 
 $H^i(\mathcal{C}(\mathcal{\MM})^\bullet)$ for $i\neq 0$ are supported
 on the irrelevant ideal and hence that 
$$\det(\mathcal{C}(\mathcal{\MM})^\bullet \otimes \mathbf{C}(H)
)=\sum_{D} \mathrm{mult}_{(Z_D)}H^0(q_*(\LL^\bullet)).D.$$
Using again that the morphism  $q$
is  birational from $W$ to the resultant divisor  
$q(W)=\nabla$, we deduce the theorem.
\end{preuve}

\begin{rem} {\rm The theorem \ref{detresdet} generalizes the construction
    given for the resultant of a very ample vector bundle  in
    \cite{GKZ94}. Indeed the complex that we obtain here 
    degenerates in the Koszul complex 
    of the universal section of the vector bundle 
    $\hom(E,\OO_X)$, by setting $F=\OO_X$ and $r=0$.
    }
\end{rem}

This result shows that the computation of the determinantal resultant
can be done by computing the determinant of a certain complex (which
can be done by using the method of Cayley, see \cite{GKZ94}, appendix B). As a
consequence the determinantal resultant is also obtained as the gcd
of all the determinants of the maximal minors of the surjective map
$$\mathcal{C}(\mathcal{\MM})^{-1}\otimes \mathbf{C}(H) \rightarrow 
\mathcal{C}(\mathcal{\MM})^{0}\otimes \mathbf{C}(H),$$
coming from the complex $\mathcal{C}(\mathcal{\MM})^\bullet \otimes
\mathbf{C}(H)$. In this way, being given $\varphi \in H$, the problem
of testing if $X_r(\varphi)$ is empty or not is traduced in a problem
of linear algebra since it corresponds to test if a matrix is
surjective or not (which is basically done by rank computations). In
the following section we will detail these results in the case 
where  $X$ is a projective space.


\section{Determinantal resultant on projective spaces}\label{Pdetres}

In this section we focus on the particular case where  
$X=\PP^{(m-r)(n-r)-1}$, where $m,n$ 
and $r$ are three positive integers such that $m\geq n > r \geq 0$,
and  where the vector bundles $E$ and $F$ are given by 
$$E=\bigoplus_{i=1}^m \OO_X(-d_i), \ F=\bigoplus_{j=1}^n\OO_X(-k_j),$$
the integers  $d_i$ and $k_j$ (not necessary positives) satisfying $d_i
> k_j$ for all  $i,j$.
The vector space $H=\hom(E,F)$ is identified with the vector space of 
matrices of size $n\times m$ with entry $i,j$  a  
homogeneous polynomial on $X$ of degree $d_j-k_i$, i.e. the 
matrices
$$\left( \begin{array}{cccc}
    h_{1,1} & h_{1,2} & \ldots & h_{1,m} \\
    h_{2,1} & h_{2,2} & \ldots & h_{2,m} \\
    \vdots & \vdots & & \vdots \\
    h_{n,1} & h_{n,2} & \ldots & h_{n,m} \\
  \end{array} \right),$$
where $h_{i,j} \in H^0(X,\OO_X(d_j-k_i))$. By theorem 
\ref{defresdet} the $\rth$-determinantal resultant 
$\Res_{E,F,r}$ exists and vanishes for all matrix with rank  lower
or equal to $r$ in at least on point of $X$. It is a generalization of
the classical Macaulay resultant (corresponding to the case $n=1$, $r=0$ and
$k_1=0$).

\begin{rem} {\rm Notice that for all $l\in \ZZ$ both determinantal
    resultants $\Res_{E,F,r}$ and
    $\Res_{E\otimes\OO(l),F\otimes\OO(l),r}$ are equals.  Consequently
    we can always suppose that one of the integers 
        $d_i$ or one of the integers $k_j$ is zero.}
 \end{rem}

In what follows we will first look at the degree of the determinantal
resultant in this case, then we will explicit its
computation and end with an example.

\subsection{The degree}
Observing proposition \ref{multideg} it appears that the degree of
the $\rth$-determinantal resultant (when it exists) is given by a
closed formula which only depends on the integers $d_i$ and
$k_j$. Indeed the Chern polynomials of $E$ and $F$ are given by the 
formulas~: 
$$c_t(E)=\prod_{i=1}^m (1-d_it) \ \ , \  c_t(F)=\prod_{i=1}^n
(1-k_it).$$
To be more precise, the degree of the $\rth$-determinantal
resultant in the coefficients of the column number $i$ (that is in the
coefficients of the polynomials $h_{1,i},h_{2,i},\ldots,h_{n,i})$ is the
coefficient of $\alpha_i$ of the multivariate polynomial (in variables
$\alpha_1,\ldots,\alpha_m$)
\begin{equation}\label{degr}
(-1)^{(m-r)(n-r)}\Delta_{m-r,n-r}\left(\frac{\prod_{i=1}^{m}(1-(d_i+\alpha_i)t)}{\prod_{i=1}^{n}(1-k_it)}\right).
\end{equation}
For instance the degree of the principal determinantal resultant is
obtained as the 
coefficient of $\alpha_i$ of the multivariate polynomial (in variables
$\alpha_1,\ldots,\alpha_m$) computed as the coefficient of the
monomial $t^{m-n+1}$
in the univariate polynomial (in the variable $t$) 
$$(-1)^{(m-r)(n-r)}\frac{\prod_{i=1}^{m}(1-(d_i+\alpha_i)t)}{\prod_{i=1}^{n}(1-k_it)}.$$
Let us see what we obtain if we take the simple example of
$X=\PP^1$. The arithmetic conditions on $m$, $n$ and $r$ implies that
$r=n-1$ (principal case) and $m-n=1$. Hence the simplest determinantal
resultant is obtained with $E=\OO(-d_1)\oplus\OO(-d_2)$ and
$F=\OO$. It is the well-known Sylvester resultant of two
polynomials. We see easily that the 
coefficient of $t^2$ in \eqref{degr} is
$(d_1+\alpha_1)(d_2+\alpha_2)$ and hence recover that the degree of
the Sylvester resultant is $d_1$ in the coefficients of the second
column and $d_2$ in the coefficients of the first one.\\
We can also look at the principal determinantal resultant
corresponding to  $E=\OO(-d_1)\oplus\OO(-d_2)\oplus\OO(-d_3)$ and
$F=\OO(-k)\oplus\OO$. Denoting by $N_i$ the degree of this resultant
in the coefficients of the column number $i$, for $i=1,2,3$, and
applying \eqref{degr} we obtain $N_1=d_2+d_3-k$, $N_2=d_1+d_3-k$ and
$N_3=d_1+d_2-k$.\\ 
Finally we  finish with an example where we know by advance the
multi-degrees. This example corresponds to the case
$E=\OO(-d_1)\oplus\OO(-d_2)$, $F=\OO(-k)\oplus\OO$ with $r=0$ (it is
not a principal case). We are
hence over $X=\PP^3$. Applying our results we can find
$N_1=d_2(d_2-k)(2d_1-k)$ and $N_2=d_1(d_1-k)(2d_2-k)$. In fact this
determinantal condition  corresponds to the classical resultant of the
polynomials $h_{1,1}, h_{1,2}, h_{2,1}, h_{2,2}$ of respective degree
$d_1-k, d_2-k, d_1$ and $d_2$ since we have $r=0$. For such a
resultant we know that the degree in the coefficient of one of these
polynomial is the product of the degrees of the others and we can in
this way check our formulas for $N_1$ and $N_2$.

\subsection{The computation}
By theorem \ref{detresdet} determinantal resultants on projective
spaces can be
 computed as the determinant of a certain complex defined in the
 preceding section. This complex involves a line bundle $\MM$ that we
 are going to precise here. In fact we have to choose this line
 bundle such that 
 $$H^i(X,\wedge_I(E)\otimes S_{I'}(F^*) \otimes \MM)=0,$$
 for all  $i>0$, and for all partition $I$ such that
 $-|I|+n(I)=p$. The integers 
 $d_j-k_i$ being all positives, we deduce that we have only to have
 that $\MM$ stabilizes the term on the far left of the complex 
 $\mathcal{L}^\bullet$ to stabilize all the complex (this is a direct
 consequence of the cohomology of line bundles on projective spaces,
 see \cite{Hart77}, theorem  III, 5.1). In this way, we have to
 determine all the integers $d$ such that 
 $$H^i(X,\wedge_I(E)\otimes S_{I'}(F^*) \otimes \OO(d))=0,$$
 for all $i>0$, and where $I$ is the partition
 $I=(\underbrace{m,\ldots,m}_{q})$ (recall that $q=n-r$).

 To do this we have to state another property of the Schur functors,
 the additivity property~: if $A$ and $B$ are two vector bundles on a
 scheme $Z$, and $I$ a partition, then 
 $$S_I(A\oplus B)\simeq \bigoplus_J S_{I/J}(A)\oplus S_J(B),$$
 where the sum is limited to partitions $J$ contained in $I$ (see
 \cite{Las78}). If now 
  $L$ is a line bundle on $Z$, we deduce that 
 for all pair of partitions ($I,H$),
 $$S_{I/H}(A\oplus L)\simeq \bigoplus_J L^{\otimes|J/H|}\otimes S_{I/J}(A),$$
 the sum being on all partitions $J$ such that 
 $$H \subset J \subset I, \ {\rm if} \ h_1\leq j_1\leq h_2 \leq j_2
 \ldots.$$
 We hence deduce, for $I=(\underbrace{m,\ldots,m}_{q})$, that 
 $$\wedge_I(E)\simeq \OO(-q(\sum_{i=1}^m d_i))
 \bigoplus_{s}\OO(-r_s),$$
 where the integers $r_s$ satisfy $r_s\leq q(\sum_{i=1}^m d_i)$.\\
 The partition $I'$ associated to the partition $I$ is defined by 
 $$I'=(\underbrace{n-r}_{r},\underbrace{m-r}_{q}).$$
 Supposing $k_1\geq k_2 \geq \ldots \geq k_n$, it follows 
 $$S_{I'}(F^*)\simeq \OO((n-r)(k_1+\ldots+k_r)+(m-r)(k_{r+1}+\ldots +k_n))
 \bigoplus_{s}\OO(t_s),$$
 where the integers $t_s$ satisfy $t_s \leq
 (n-r)(k_1+\ldots+k_r)+(m-r)(k_{r+1}+\ldots +k_n)$.\\
 The cohomology of line bundles on projective spaces 
 shows that the integer $d$ must be such that~:
 $$d - q(\sum_{i=1}^m d_i) + (n-r) (\sum_{i=1}^n k_i) +
 (m-n)(k_{r+1}+\ldots k_n) \geq -(m-r)(n-r)+1,$$
 that is 
 $$d \geq q(\sum_{i=1}^m d_i) - (n-r) (\sum_{i=1}^n k_i) -
 (m-n)(k_{r+1}+\ldots+k_n) - (m-r)(n-r)+1.$$
 We denote by $\nu_{\dg,\kg}$ the integer on the right of this last
 equality, that is 
 $$\nu_{\dg,\kg}=(n-r)(\sum_{i=1}^m d_i - \sum_{i=1}^n k_i) -
 (m-n)(k_{r+1}+\ldots+k_n) - (m-r)(n-r)+1.$$
\begin{prop} Suppose that $k_1\geq k_2 \geq \ldots \geq k_n$, then 
for all integer $d\geq \nu_{\dg,\kg}$, we have
$$\det(\mathcal{C}(\OO_{X}(d))^\bullet \otimes \mathbf{C}(H) 
)=\Res_{E,F,r}.$$ 
\end{prop}

 \begin{rem} {\rm 
      The determinantal resultant is invariant if we twist 
        both vector bundles  $E$ and $F$ by the same line bundle 
        $\OO(l)$, with $l\in \ZZ$. We can check by an easy
        computation that the integer $\nu_{\dg,\kg}$  we just
        defined is invariant by this transformation, i.e. 
        $$\nu_{\dg,\kg}(E,F)=\nu_{\dg,\kg}(E\otimes \OO(l),F\otimes
        \OO(l)).$$
        Notice also that the formula of the integer $\nu_{\dg,\kg}$
        for the determinantal resultant is 
        $$\nu_{\dg,\kg}=(d_1+\ldots+d_m - m) - (k_1+\ldots+k_n - n).$$
        Here again, for the case $n=1$ we recover that
        $\nu_{\dg,\kg}$ is the known critical degree for the Macaulay 
        resultant, that is $d_1+\ldots+d_m -m +1$ (see \cite{Jou97}),
        since we can suppose that $k_1=k_n=0$ without changing the
        determinantal resultant. 
}
 \end{rem}

We can now explicit how to compute the determinantal resultant.  
Denoting $R=\KK[x_0,\ldots,x_{(m-r)(n-r)-1}]$, the first map 
$\sigma_{d}$ (the one of the far right) of the complex
$\mathcal{C}(\OO_{X}(d))^\bullet \otimes \mathbf{C}(H)$ is the map  
$$\begin{array}{ccl}
  \bigoplus_{i_1<\ldots < i_{r+1},\, j_1<\ldots <j_{r+1}}
  R_{[d-\sum_{t=1}^{r+1}d_{i_t}+\sum_{t=1}^{r+1}k_{i_t}]}e_{i_1,\ldots ,i_{r+1},j_1,\ldots ,j_{r+1}} &
  \xrightarrow{\sigma_{d}} & R_{[d]} 
\end{array}$$
which associates to each $e_{i_1,\ldots, i_{r+1},j_1,\ldots ,j_{r+1}}$
the polynomial $\Delta_{i_1,\ldots ,i_{r+1},\, j_1,\ldots ,j_{r+1}}$ 
 denoting the determinant of the minor 
$$\left( \begin{array}{cccc}
    h_{j_1,i_1} & h_{j_1,i_2} & \ldots & h_{j_1,i_{r+1}} \\
    h_{j_2,i_1} & h_{j_2,i_2} & \ldots & h_{j_2,i_{r+1}} \\
    \vdots & \vdots & & \vdots \\
    h_{j_{r+1},i_1} & h_{j_{r+1},i_2} & \ldots & h_{j_{r+1},i_{r+1}} \\
  \end{array} \right),$$
$R_{[t]}$ being the vector space of homogeneous polynomials of fixed
degree $t$. By properties of complex determinants 
 (see \cite{GKZ94}, appendix A), we deduce the following algorithm~:

\begin{prop} Choose an integer $d\geq \nu_{\dg,\kg}$. All nonzero maximal
  minor (of size $\sharp R_{[d]}$) of the map $\sigma_{d}$ is a  
 multiple of the $\rth$-determinantal resultant $\Res_{E,F,r}$. 
  Moreover the greatest common divisor of all the determinants of
  these maximal minors is exactly $\Res_{E,F,r}$.
\end{prop}

This proposition gives us an algorithm to compute explicitly 
the determinantal resultant, completely similar  
to the one giving the expression of the Macaulay resultant. Notice
that it is also possible to give the equivalent (in a less 
explicit form) of the 
so-called Macaulay matrices (of the Macaulay resultant) for the
principal determinantal resultant. These Macaulay matrices are in fact some
particular maximal minors of the map $\sigma_d$ of the
Macaulay resultant, i.e  $n=1$ and $r=0$. They are obtained by
specializing the  input polynomials in $x_0^{d_0}, x_1^{d_1}, \ldots,
x_m^{d_m}$ (see \cite{Jou97}). For the principal determinantal
resultant ($r=n-1$) it is also possible to give a specialization of polynomials
$(h _{i,j})$ so that we can obtain maximal minors explicitly. The specialization is the 
following one~:
$$\left(\begin{array}{cccccccc}
  x_0^{d_1-k_1} & x_1^{d_2-k_1} & \cdots & x_{m-n}^{d_{m-n+1}-k_1} & 0 &
  \cdots & \cdots & 0 \\ 
  0 & x_0^{d_2-k_2} & x_1^{d_3-k_2} & \cdots & x_{m-n}^{d_{m-n+2}-k_2}
  & 0 & \cdots & 0 \\ 
  \vdots & \vdots & \vdots & \vdots & \vdots & \vdots & \vdots & \vdots \\
  0 & \cdots & 0 & x_0^{d_n-k_n} & x_1^{d_{n+1}-k_n} & \cdots & \cdots 
  & x_{m-n}^{d_m-k_n}  \\
\end{array}\right).$$
It is easy to check that the rank of this matrix is $n$ at all points
of $X=\PP^{m-n}$ and hence that its maximal minors can be left to the
matrix $\sigma_d$. 

However, from a computational point of view we have to notice that developed
determinants are not efficient. In practice it is preferable to work
with the matrix of $\sigma_d$ itself. For instance to see if some
particular values of the parameters vanish the determinantal resultant we just
have to put them in the matrix and check its rank.

\subsection{Chow forms of rational normal scrolls}
We now illustrate the determinantal resultant by showing how it
computes the Chow form of a rational normal scroll. Notice that such
formulas are exposed in \cite{EiSc01}, examples 2.4 and 2.5.

A rational normal scroll can be basically defined by its equations. Let
$d_1,\ldots,d_r$ be $r$ positive integers such that at least one of
them is strictly positive. Denoting $N+1=\sum_{i=1}^{r}(d_i+1)$, we
take the homogeneous coordinates on $\PP^N$ to be
$$X_{1,0},X_{1,1},\ldots,X_{1,d_1}, X_{2,0},\ldots,X_{2,d_2},
\ldots, X_{r,0},\ldots,X_{r,d_r}.$$
Define a $2\times(\sum_{i=1}^{r}d_i)$ matrix $M(d_1,\ldots,d_r)$ of
linear forms on 
$\PP^N$ by
$$\left(\begin{array}{cccccccccc}
X_{1,0} & \ldots & X_{1,d_1-1} & X_{2,0} & \ldots & X_{2,d_2-1} &
\ldots & X_{r,0} & \ldots & X_{r,d_r-1} \\
X_{1,1} & \ldots & X_{1,d_1} & X_{2,1} & \ldots & X_{2,d_2} &
\ldots & X_{r,1} & \ldots & X_{r,d_r}
\end{array}\right).$$
The rational normal scroll $S(d_1,\ldots,d_r) \subset \PP^N$ is the
variety defined by the ideal of $2\times 2$ minors of
$M(d_1,\ldots,d_r)$. It is an irreducible variety of dimension $r$ and
degree $\sum_{i=1}^{r}d_i$, as it is proved in \cite{EiHa87} (notice
that $\deg(S)=1+\codim(S)$, a rational normal scroll is a variety of
minimal degree). In this paper another (equivalent) definition is
given. Let $F$ be the vector bundle
$F=\oplus_{i=1}^{r}\OO_{\PP^1}(d_i)$. By
hypothesis we made on the integers $d_i$, $i=1,\ldots,r$, the vector
bundle $F$ is generated by $N+1$ global sections. The
projectivized vector bundle $\PP(F)$ is a smooth variety of
dimension $r$ mapping to $\PP^1$ with fibers $\PP^{r-1}$. Its
tautological line bundle $\OO_{\PP(F)}(1)$ is generated by its global
sections and defines a tautological map $\PP(F)\rightarrow
\PP^N$. This map is birational and $S(d_1,\ldots,d_r)$ is its image.

Recall now that the Chow divisor of a $k$-dimensional variety
$X\subset \PP^n$ is the hypersurface in the Grassmannian $G$ of planes
of codimension $k+1$ in $\PP^n$, consisting of those planes meeting
$X$. The Chow form of $X$, denoted $\mathcal{C}(X)$, is its defining
equation, defined up to the multiplication by a non-zero constant. If
$X$ is of degree $e$ then $\mathcal{C}(X)$ is also of degree $e$ in
$G$. From the definition of the rational normal scroll $S(d_1,\ldots,d_r)$
we deduce that its Chow form $\mathcal{C}(S(d_1,\ldots,d_r)))$ can be
computed as a principal determinantal resultant (for which we do not
have to suppose that th ground field is of characteristic zero). Indeed, a
$(r+1)$-dimensional space of sections of $F$, say $\alpha : \OO^{r+1}
\rightarrow F$, corresponds to a plane of codimension $r+1$ in $\PP^N$
which meets $X$ if and only if there exists a point $p$ of $\PP^1$
such that $\rank(\alpha(p))\leq r$. We have the following result~:

\begin{prop}\label{chow}
The Chow form of the rational normal scroll $S(d_1,\ldots,d_r)$ is
given by the principal determinantal resultant $\Res_{E,F,r}$, where
$E=\OO^{r+1}_{\PP^1}$ and $F=\oplus_{i=1}^{r}\OO_{\PP^1}(d_i)$.
\end{prop}

Let us illustrate this proposition with the example of the rational
normal scroll $S(2,1)$ (see also \cite{EiSc01}, example 2.5). The
``generic'' map $\alpha : \OO^3 \rightarrow \OO(2)\oplus \OO(1)$ is
given by the matrix
$$\begin{pmatrix}{{a}}_{0} x^{2}+{{a}}_{1} {x} {y}+{{a}}_{{2}} y^{2}&
     {{b}}_{0} x^{2}+{{b}}_{1} {x} {y}+{{b}}_{{2}} y^{2}&
     {{c}}_{0} x^{2}+{{c}}_{1} {x} {y}+{{c}}_{{2}} y^{2}\\
     {{a}}_{{3}} x+{{a}}_{{4}} y&
     {{b}}_{{3}} x+{{b}}_{{4}} y&
     {{c}}_{{3}} x+{{c}}_{{4}} y\\
     \end{pmatrix},$$
where $(x:y)$ denotes the homogeneous coordinates of $\PP^1$.
The matrix associated to its principal determinantal resultant ($r=1$) is given by the $5\times
6$ matrix~:
{\scriptsize 
$$\left(\begin{array}{ccc} {{a}}_{{3}} {{b}}_{0}-{{a}}_{0} {{b}}_{{3}}&
      0&
      {{a}}_{{3}} {{c}}_{0}-{{a}}_{0} {{c}}_{{3}}\\
      {{a}}_{{4}} {{b}}_{0}+{{a}}_{{3}} {{b}}_{1}-{{a}}_{1} {{b}}_{{3}}-{{a}}_{0} {{b}}_{{4}}&
      {{a}}_{{3}} {{b}}_{0}-{{a}}_{0} {{b}}_{{3}}&
      {{a}}_{{4}} {{c}}_{0}+{{a}}_{{3}} {{c}}_{1}-{{a}}_{1} {{c}}_{{3}}-{{a}}_{0} {{c}}_{{4}}\\
      {{a}}_{{4}} {{b}}_{1}+{{a}}_{{3}} {{b}}_{{2}}-{{a}}_{{2}} {{b}}_{{3}}-{{a}}_{1} {{b}}_{{4}}&
      {{a}}_{{4}} {{b}}_{0}+{{a}}_{{3}} {{b}}_{1}-{{a}}_{1} {{b}}_{{3}}-{{a}}_{0} {{b}}_{{4}}&
      {{a}}_{{4}} {{c}}_{1}+{{a}}_{{3}} {{c}}_{{2}}-{{a}}_{{2}} {{c}}_{{3}}-{{a}}_{1} {{c}}_{{4}}\\
      {{a}}_{{4}} {{b}}_{{2}}-{{a}}_{{2}} {{b}}_{{4}}&
      {{a}}_{{4}} {{b}}_{1}+{{a}}_{{3}} {{b}}_{{2}}-{{a}}_{{2}} {{b}}_{{3}}-{{a}}_{1} {{b}}_{{4}}&
      {{a}}_{{4}} {{c}}_{{2}}-{{a}}_{{2}} {{c}}_{{4}}\\
      0&
      {{a}}_{{4}} {{b}}_{{2}}-{{a}}_{{2}} {{b}}_{{4}}&
      0\\
      \end{array}\right.$$

$$\left.\begin{array}{ccc} 0&
      {{b}}_{{3}} {{c}}_{0}-{{b}}_{0} {{c}}_{{3}}&
      0\\
      {{a}}_{{3}} {{c}}_{0}-{{a}}_{0} {{c}}_{{3}}&
      {{b}}_{{4}} {{c}}_{0}+{{b}}_{{3}} {{c}}_{1}-{{b}}_{1} {{c}}_{{3}}-{{b}}_{0} {{c}}_{{4}}&
      {{b}}_{{3}} {{c}}_{0}-{{b}}_{0} {{c}}_{{3}}\\
      {{a}}_{{4}} {{c}}_{0}+{{a}}_{{3}} {{c}}_{1}-{{a}}_{1} {{c}}_{{3}}-{{a}}_{0} {{c}}_{{4}}&
      {{b}}_{{4}} {{c}}_{1}+{{b}}_{{3}} {{c}}_{{2}}-{{b}}_{{2}} {{c}}_{{3}}-{{b}}_{1} {{c}}_{{4}}&
      {{b}}_{{4}} {{c}}_{0}+{{b}}_{{3}} {{c}}_{1}-{{b}}_{1} {{c}}_{{3}}-{{b}}_{0} {{c}}_{{4}}\\
      {{a}}_{{4}} {{c}}_{1}+{{a}}_{{3}} {{c}}_{{2}}-{{a}}_{{2}} {{c}}_{{3}}-{{a}}_{1} {{c}}_{{4}}&
      {{b}}_{{4}} {{c}}_{{2}}-{{b}}_{{2}} {{c}}_{{4}}&
      {{b}}_{{4}} {{c}}_{1}+{{b}}_{{3}} {{c}}_{{2}}-{{b}}_{{2}} {{c}}_{{3}}-{{b}}_{1} {{c}}_{{4}}\\
      {{a}}_{{4}} {{c}}_{{2}}-{{a}}_{{2}} {{c}}_{{4}}&
      0&
      {{b}}_{{4}} {{c}}_{{2}}-{{b}}_{{2}} {{c}}_{{4}}\\
      \end{array}\right).$$}

We can compute the gcd of the maximal minors (size $5\times 5$) of
this matrix and we obtain the Chow form of $S(2,1)$ in the Stiefel
coordinates (see \cite{GKZ94}, chapter 3.1 for different coordinates
of Grasmannians). It is of degree $3$ in each set of variables
$(a_i)_{i=0\ldots 4}$, $(b_i)_{i=0\ldots 4}$ and  $(c_i)_{i=0\ldots
  4}$~:
$$-a_4^3b_2^2b_3c_0^3+a_3a_4^2b_2^2b_4c_0^3+2a_2a_4^2b_2b_3b_4c_0^3
-2a_2a_3a_4b_2b_4c_0^3-a_2^2a_4b_3b_4^2c_0^3+\ldots \ .$$
The Pl\"ucker coordinates are obtained from the Stiefel
coordinates as the $3\times 3$ minors of
the matrix
$$\left(\begin{array}{ccccc}
    a_0 & a_1 & a_2 & a_3 & a_4 \\
    b_0 & b_1 & b_2 & b_3 & b_4 \\
    c_0 & c_1 & c_2 & c_3 & c_4
\end{array}\right),$$
and we can check that the Chow form of $S(2,1)$ is well of total
degree 3 in these coordinates. From a computational point of view it is
not efficient to compute such an extended formula for
$\mathcal{C}(S(2,1))$. It is better to work with a matricial representation of
this Chow form, here with the $5\times 6$ matrix we have computed, and
use it. For instance, to test if a point of the grassmannian (in the
Stiefel coordinates) is on $\mathcal{C}(S(2,1))$ is  
done by replacing the $a_i$'s $b_i$'s and $c_i$'s by their value and
 computing the rank of the matrix~: if the rank is 5 the point is
not on the Chow form, otherwise it is.

Looking to the general case of proposition \ref{chow}, the Chow form of
the rational normal scroll $S(d_1,\ldots,d_r)$ is obtained in this way
in the Stiefel coordinates of the Grassmannian of planes of
codimension $r+1$ in $\PP^N$. By \eqref{degr} we see that the
principal determinant we compute is of degree $\sum_{i=1}^rd_i$ with
respect to each column of the ``generic'' map $\OO_{\PP^1}^{r+1}\rightarrow
\oplus_{i=1}^r \OO_{\PP^1}(d_i)$, and as a polynomial of total degree
$\sum_{i=1}^rd_i$ (the degree of $S(d_1,\ldots,d_r)$) in the Pl\"ucker
coordinates.


\end{document}